\documentclass[12pt, dvips, twoside]{amsart}
\usepackage{amsmath,amssymb,amscd,amsxtra,latexsym,%bbm,
graphicx,epsfig,%epic,eepic,oldgerm,bm,
euscript,%psfrag,
amsthm}
\usepackage{a4wide}
\usepackage{amsmath,amsfonts,amsthm,amsopn,amssymb,enumitem}
\usepackage{palatino}
\usepackage[dvips]{color}
\theoremstyle{plain}
\newtheorem{theorem}{Theorem}[section]
\newtheorem*{theorem*}{Theorem}
\newtheorem{lemma}[theorem]{Lemma}
\newtheorem{proposition}[theorem]{Proposition}

\theoremstyle{definition}
\newtheorem*{remarks*}{Remark}
\newtheorem{remark}[theorem]{Remark}

\newtheorem*{example*}{Example}
\newtheorem*{examples*}{Examples}
\newtheorem{definition}[theorem]{Definition}
\newtheorem*{definition*}{Definition}

\newcommand{\proofend}{\hspace*{\fill} $\Box$\\}

\def\1{\:\!}
\def\2{\;\!}

\def\Diffc0{\operatorname{Diff^c_0}}

\def\Sympc0{\operatorname{Symp^c_0}}

\def\pp{\partial}

\def\ni{\noindent}
\def\b{\bigskip}

\def\.{\mskip1mu}
\def\?{\mskip-1mu}

\def\1{\:\!}
\def\2{\;\!}

\def\Diffc0{\operatorname{Diff^c_0}}

\def\Sympc0{\operatorname{Symp^c_0}}

\def\pp{\partial}

\def\ni{\noindent}
\def\b{\bigskip}

\def\.{\mskip1mu}
\def\?{\mskip-1mu}

\def\R{\mathbb{R}}

\def\tf{\tilde{f}}
\def\tF{\tilde{F}}

\begin{document}

\title{A rigidity result for overdetermined elliptic problems in the plane}

\author{Antonio Ros}
\address{(A.~Ros)
Departamento de Geometr\'ia y Topolog\'ia,
Universidad de Granada,
Campus Fuentenueva,
18071 Granada,
Spain}
\email{aros@ugr.es}

\author{David Ruiz}
\address{(D.~Ruiz)
Departamento de An\'alisis matem\'atico, Universidad de Granada,
Campus Fuentenueva, 18071 Granada, Spain} \email{daruiz@ugr.es}

\author{Pieralberto Sicbaldi}
\address{(P.~Sicbaldi)
Aix Marseille Universit\'e - CNRS, Centrale Marseille - I2M, UMR 7373, 13453 Marseille, France}
\email{pieralberto.sicbaldi@univ-amu.fr}

\date{\today}

\maketitle
\begin{abstract}
Let $f:[0,+\infty) \to \mathbb{R}$ be a (locally) Lipschitz
function and $\Omega \subset \mathbb{R}^2$ a $C^{1,\alpha}$ domain
whose boundary is unbounded and connected. If there exists a
positive bounded solution to the overdetermined elliptic problem
\[
\left\{\begin{array} {ll}
\Delta u + f(u) = 0 & \mbox{in }\; \Omega\\
               u= 0 & \mbox{on }\; \pp \Omega \\
\frac{\partial u}{\partial \vec{\nu}}=1 &\mbox{on }\; \pp \Omega
\end{array}\right.
\]
we prove that $\Omega$ is a half-plane. In particular, we obtain a
partial answer to a question raised by H. Berestycki, L. Caffarelli
and L. Nirenberg in 1997.
\end{abstract}

\section{Introduction}  \label{s:intro}

Given a locally Lipschitz function $f$, a widely open problem is
to classify the set of domains $\Omega \subset \R^n$ where there
exists a bounded solution $u$ to the overdetermined elliptic
problem
\begin{equation}\label{pr_bis}
\left\{\begin{array} {ll}
\Delta u + f(u) = 0 & \mbox{in }\; \Omega\\
                        u > 0 & \mbox{in }\; \Omega\\
               u= 0 & \mbox{on }\; \pp \Omega \\
\frac{\partial u}{\partial \vec{\nu}}=1 &\mbox{on }\; \pp
\Omega\,.
\end{array}\right.
\end{equation}
Here $\vec{\nu}(x)$ stands for the interior normal vector to
$\partial \Omega$ at $x$. In this case we say that $\Omega$ is an
{\it $f$-extremal domain} (see \cite{Ros-Sic} for a motivation of
that definition). The case of bounded $f$-extremal domains was
completely solved by J. Serrin in \cite{Serrin} (see also
\cite{Puc-Ser}): the ball is the unique such domain and any
solution is radial. This result has many applications to Physics
and Applied Mathematics (see \cite{Far-Val-1, Fre, Ros-Sic, Sok}).
Instead, the case of unbounded domains $\Omega$ is far from being
completely understood.

\medskip

Overdetermined boundary conditions arise naturally in free
boundary problems,  when the variational structure imposes
suitable conditions on the separation interface, see for example
\cite{Alt-Caf}. In this context, several methods for studying the
regularity of the interface are based on blow-up techniques which
lead to the study of an elliptic problem in an unbounded domain.
In this framework, problem (\ref{pr_bis}) in unbounded domains was
considered in \cite{BCN} for $f(u) = u-u^3$ (the Allen-Cahn
equation). In that paper, H. Berestycki, L. Caffarelli and L.
Nirenberg proposed the following:

\medskip

{\bf Conjecture BCN}: If $\mathbb{R}^n \backslash
\overline{\Omega}$ is connected, then the existence of a bounded
solution to problem \eqref{pr_bis} implies that $\Omega$ is either
a ball, a half-space, a generalized cylinder $B^k \times
\mathbb{R}^{n-k}$ ($B^k$ is a ball in $\mathbb{R}^k$), or the
complement of one of them.

\medskip

That question was motivated by the results of the same authors in
\cite{BCN}, and some other results concerning exterior domains,
i.e. domains that are the complement of a compact region (see
\cite{aft-bus, reichel}).

\medskip

In \cite{Sicbaldi} the third author gave a counterexample to that
conjecture for $n \geq 3$, constructing a periodic perturbation of
the straight cylinder $B^{n-1} \times \mathbb{R}$ that supports a
periodic solution to problem (\ref{pr_bis}) with $f(t) = \lambda\,
t$. The goal of this paper is to prove that Conjecture BCN is true
for $n=2$ if $\partial \Omega$ is unbounded.

\medskip

In the last years, a deep parallelism between overdetermined
elliptic problems and constant mean curvature (CMC) surfaces has
been observed. Serrin's result can be seen as the analogue of the
Alexandrov's one (\cite{Alex}), which asserts that the only
embedded compact CMC hypersurfaces in $\R^{n}$ are round spheres.
In \cite{Sch-Sic} F. Schlenk and the third author show that the
counterexamples to Conjecture BCN built in \cite{Sicbaldi} belong
to a smooth one-parameter family that can be seen as a counterpart
of the family of {\it Delaunay surfaces}. In \cite{traizet} M.
Traizet finds a one-to-one correspondence between $0$-extremal
domains in dimension $2$ and a special class of minimal surfaces
(see Section 2 for the exact statement of this result). In
\cite{DPW} M. Del Pino, F. Pacard and J. Wei consider problem
\eqref{pr_bis} for functions $f$ of Allen-Cahn type and they build
new solutions in domains in $\R^3$ whose boundary is close to a
dilated Delaunay surface or a dilated minimal catenoid. They also
build bounded and monotone solutions to problem (1) for epigraphs
in case $n \geq 9$ (this type of solutions do not exist if $n \leq
8$, as has been proved by K. Wang and J. Wei in \cite{WW}). The
domain in \cite{DPW} has boundary close to a dilated Bombieri-De
Giorgi-Giusti entire minimal graph (\cite{BDG}).

\medskip

We point out that almost all those examples of $f$-extremal
domains have boundary with some nontrivial topology. The only
exception is the epigraph extremal domain found in \cite{DPW},
which requires $n\geq 9$. Therefore it is natural to consider BCN
Conjecture if $\partial \Omega$ has the topology of the Euclidean
space and $n \leq 8$. In this paper we solve the case $n=2$.

\medskip

Some partial results have been already given in the literature for
dimension $2$. In \cite{Far-Val-0} A. Farina and E. Valdinoci
prove BCN Conjecture if $u$ is monotone along one direction and
$\nabla u$ is bounded. In \cite{WW} the case of $f$-extremal
epigraphs is solved for some nonlinearities $f$ of the Allen-Cahn
type. Finally, in \cite{Ros-Sic} the result is proved if either
$f(t) \geq \lambda t$ or $\Omega$ is contained in a half-plane and
$\nabla u$ is bounded (see also \cite{espinar} for a
generalization to other geometries). Observe that the assumption
$f(t) \geq t$ excludes the prototypical Allen-Cahn nonlinearity;
we point out that not even the half-plane is an $f$-extremal
domain for those nonlinearities $f$. In this paper we prove
Conjecture BCN for $n=2$ under the only assumption that $\partial
\Omega$ is unbounded. The exact statement of our result is the
following:

\medskip

\begin{theorem} \label{Tmain} Let $\Omega \subset \mathbb{R}^2$ a $C^{1,\alpha}$
domain whose boundary is unbounded and connected. Assume that
there exists a bounded solution $u$ to problem \eqref{pr_bis} for
some (locally) Lipschitz function $f:[0,+\infty) \to \mathbb{R}$.
Then $\Omega$ is a half-plane and $u$ is parallel, that is, $u$
depends only on one variable.
\end{theorem}

We point out that, generally speaking, $f$-extremal domains always
have $C^{2,\alpha}$ regularity, as shown in \cite{vogel}. Hence,
Theorem \ref{Tmain} could be stated under less regularity
requirements, but for the sake of clarity we have preferred to
leave it in that form.

\medskip

The proof is divided in several steps. First, we show that the
curvature of $\partial \Omega$ is bounded. This is done via a
blow-up argument, making use of the classification results of
\cite{traizet} for the case $f=0$. This argument needs some
uniform regularity estimates that are given in Section 2, together
with other preliminary results. In particular this result implies,
via standard regularity for elliptic problems, that $\nabla u$ is
bounded. This allows us to prove Theorem \ref{Tmain} if $u$ is
monotone along one direction. This result is basically contained
in \cite{Far-Val-0} if $\partial \Omega$ is $C^3$; in Section 4 we
relax this regularity assumption by using ideas from the proof of
the De Giorgi conjecture in dimension 2. In Section 5 we combine
the previous result and the moving plane method (as well as the
so-called {\it tilted moving plane method}) to show that $\Omega$
must contain a half-plane. A crucial ingredient in our proof is
given in Section 6: we prove the existence of a divergent sequence
of points in $p_n \in \partial \Omega$ such that $\partial \Omega$
converges to a straight line near such sequence. In particular, a
parallel solution in a half-plane exists, which is given as the
limit of $u(\cdot- p_n)$. In Section 7 we use the variational
method to construct solutions in large balls converging to the
parallel solution as the radius goes to $+\infty$. Section 8
concludes the proof of Theorem \ref{Tmain}. First we show that the
graph of $u$ is above the graphs of those solutions defined in
balls: passing to the limit, it is above the parallel solution
too. But both solutions are in contact and have the same boundary
conditions, so Theorem \ref{Tmain} follows from the maximum
principle.

\b \ni {\bf Acknowledgments.} This paper was written when P. S.
visited the University of Granada during his period of
``D\'el\'egation CNRS". A. R. has been partially supported by
Mineco-Feder Grant MTM2011-22547 and by J. Andaluc{\'\i}a FQM-325. D.
R. has been supported by the Mineco Grant MTM2011-26717 and by J.
Andalucia (FQM 116). P. S. was partially supported by the GDRE
network on Geometric Analysis and ANR-11-IS01-0002 grant.

\section{Preliminary tools}

In this section we discuss some preliminary results that will be
useful throughout the paper. Throughout the paper, $B_R(p)$ stands
for the open ball of center $p$ and radius $R$.

\subsection{$C^{2,\alpha}$ regularity}

%%%xxx explicar quien seria u en este parrafo: u una overdetermined solucion de (1)
In this paper we assume that the boundary of our domains is of
class $C^{1,\alpha}$. Standard regularity arguments for elliptic
equations show that a solution $u$ of \eqref{pr_bis} is
$C^{2,\alpha}$ in $\Omega$ and $C^{1,\alpha}$ up to the boundary.
However, $f$-extremal domains always exhibit more regularity,
namely $C^{2,\alpha}$. Moreover, the following uniform estimate
holds.

\begin{lemma} \label{regularity} Fix $R>0$, $\alpha \in (0,1)$, $p=(p_1, p_2) \in \partial \Omega$ and let $\phi \in
C^{1,\alpha}(p_1-R, p_1+R)$ be such that $\Gamma_R = \partial
\Omega \cap B_R(p) \subset \{(x, \phi(x)); \ x \in (p_1-R,
p_1+R)\}$. Define $\Omega_R= \Omega \cap B_R(p)$. Let u be a
bounded solution of the problem:
\begin{equation}
\left\{\begin{array} {ll}
\Delta u =  h(x) & \mbox{in }\; \Omega\\
                        u > 0 & \mbox{in }\; \Omega\\
               u= 0 & \mbox{on }\; \pp \Omega \\
               \frac{\partial u}{\partial \vec{\nu}}=1
 &\mbox{on }\; \pp \Omega
\end{array} \right.
\end{equation}
for some $h \in C^{0,\alpha}$. Take $M= \|
h\|_{C^{0,\alpha}(\Omega_R)}+ \| u\|_{C^0(\Omega_R)} + \|
\phi\|_{C^{1,\alpha}(p_1-R, p_1+R)} $. Then, $u$, $\phi$ belong to
$C^{2,\alpha}$ and
$$ \| u\|_{C^{2,\alpha}(\Omega_{R/2})}+ \|
\phi\|_{C^{2,\alpha}(p_1-R/2, p_1+R/2)} \leq C,$$
for some $C>0$ depending only on $M$, $R$.

\end{lemma}

\begin{remark} The $C^{2,\alpha}$ regularity for overdetermined
problems in this fashion was given in \cite{KN} (see also
\cite{vogel}). However, the result in \cite{KN} needs some
additional conditions that do not hold under our assumptions.
Moreover, Lemma \ref{regularity} is also concerned with the
uniformity of the regularity estimate, which will be crucial later
on. The proof we give here is different from \cite{KN} and takes
advantage of some regularity results for problems with nonlinear
oblique boundary conditions (see \cite{lieb-book}). It is also
worth pointing out that Lemma \ref{regularity} is valid for any
dimension $n$.
\end{remark}

{\it Proof:} By standard regularity results, we conclude that
$$ \| u\|_{C^{1,\alpha}(\Omega_{2R/3})} \leq C,$$
with $C$ depending on $M$, $R$ (see \cite{gilbarg}, Theorem 8.33
and the comment that follows, and also Corollaries 8.35, 8.36).
Then, we are under the hypotheses of \cite{lieb-book}[Proposition
11.21]\footnote{In \cite{lieb-book}[Proposition 11.21] the
estimates are written with respect to a certain weighted Holder
norms, and those weights vanish when a point approaches $\partial
B_R(p) \cap \Omega$. We avoid the use of those norms by
considering estimates in a smaller ball $B_{R/2}(p)$. Moreover,
$b(x,u, \nabla u)= |\nabla u|^2-1$ and the obliquity condition
(11.57b) trivially holds in our setting.}. Therefore, there exists
$C>0$ depending on $M$, $R$ with
$$ \| u\|_{C^{2,\alpha}(\Omega_{R/2})}\leq C.$$
Now observe that $\Gamma_R$ is the $0$ level of $u$, and $|\nabla
u|=1$ there. The implicit function theorem implies that, enlarging
$C$ if necessary,
$$\|\phi\|_{C^{2,\alpha}(p_1-R/2, p_1+R/2)} \leq C.$$
This concludes the proof of the lemma. \proofend

\subsection{The moving plane method in unbounded domains}

One of the most important tools coming from the maximum principle
of elliptic operators is the {\it moving plane method}, introduced
firstly by A. D. Alexandrov \cite{Alex} for constant mean
curvature surfaces and then adapted by J. Serrin \cite{Serrin} to
elliptic overdetermined problems (see also \cite{BN, Dancer, GNN}).

\medskip

Let $L$ be a line in $\mathbb{R}^2$ that intersects $\Omega$, and
let $L^+$ and $L^-$ be the two connected components of
$\mathbb{R}^2 \backslash L$. Let us suppose that $\Omega \cap L^-$
has a bounded connected component $C$ (Fig. 1).

\begin{figure}[!ht]
\centering
{\scalebox{.5}{\input{ ./serrin6.pstex_t}}}
\label{fig1}\caption{}
\end{figure}

It is easy to prove that:
\begin{enumerate}
\item[i.] the closure of $\pp C \cap L^-$ is a graph over $\pp C
\cap L$, \item[ii.]  the closure of $\pp C \cap L^-$ is not
orthogonal to $L$ at any point, \item[iii.]  If $C'$ is the
reflection of $C$ about $L$, then the closure of $C \cup C'$ stays
within $\overline{\Omega}$, \item[iv.]  If for every $x\in C$ we
define $u'(x') = u(x)$ where $x'$ is the symmetric point to $x$
with respect to $L$, then the graph of the function $u'$ over $C'$
stays under the graph of $u$, and the two graphs are not tangent
in the points of $L$,

\item[v.]  $ \frac{\partial u}{\partial \vec{n}}>0$ in $C$ where
$\vec{n}$ is the normal direction to $L$ pointing towards $L^+$.
\end{enumerate}
The proof of these facts is a simple application of the moving plane method, and is given in \cite{Ros-Sic}.

In this case, we will say that {\it the moving plane method applies to $C$ in $L^-$  with respect to lines parallel to $L$}. We give then the following definition, where we generalize this expression to the case where $C$ is not supposed to be connected and bounded.

\begin{definition}
Let $\Omega$ be an $f$-extremal domain in $\mathbb{R}^n$. Let $L$ be a hyperplane intersecting $\Omega$, $L^-$ be one of the two components of $\mathbb{R}^n \backslash L$, and $C = \Omega \cap L^-$. We say that {\it the moving plane method applies to $C$ in $L^-$  with respect to lines parallel to $L$} when properties i.-ii.-iii.-iv.-v. above are all satisfied.
\end {definition}

An immediate consequence of the moving plane method is the
following:

\begin{lemma}\label{np}

Let $\Omega \subset \mathbb{R}^2$ be an unbounded $f$-extremal
domain such that $\partial \Omega$ is connected. Then for any
point $p \in \partial \Omega$, the half-line $N(p)$ given by the
half-line starting at $p$ and pointing in the inward normal
direction about $\partial \Omega$ with respect to $\Omega$, is
contained in $\Omega$. We say that $\Omega$ has {\it the property
of the inward normal half-line}.

\begin{figure}[!ht]
\centering
{\scalebox{.5}{\input{ ./no.pstex_t}}}
\caption{}
\label{fig7}
\end{figure}
\end{lemma}

{\it Proof}. This follows immediately from property (iii).\proofend

The property of the inward normal half-line has an interesting
consequence, that will be exploited later on:

\begin{lemma} \label{otro} Let $\Omega \subset \mathbb{R}^2$ be an unbounded $f$-extremal domain such
that $\partial \Omega$ is connected, $p \in
\partial \Omega$ and $R>0$. Denote $D$ the connected component of
$\Omega \cap B_R(p)$ with $p$ in its boundary. Then $\partial D
\cap B_R(p)$ is connected.
\end{lemma}

{\it Proof}.  Otherwise, let us call $U$ the connected component
of $\partial D \cap B_R(p)$ containing $p$. There exists $q \in
\partial D \cap B_R(p) \setminus U$ that minimizes the distance from $p$.
Clearly, the segment $[p, q]$ touches $q$ perpendicularly.

We now claim that the points of $[p,q]$ close to $q$ belong to
$\Omega$. Define $U'$ the connected component of $\partial D \cap
B_R(p)$ containing $q$. Clearly, $U'$ separates $B_R(p)$ in two
connected components, $V$ and $V'$, and $p$ belongs to one of
them, say, $V$. Since $D$ is connected, $D \subset V$. The claim
follows from the perpendicular intersection of $[p,q]$ and $U$.

By the property of the inward normal half-line, $[p,q] \subset
N(q)$ but this is a contradiction because $p \in \partial \Omega$.
%Therefore there exists $q' \in \partial D \cap N(q)$ such $[q,q'] \subset \overline{\Omega}$. This is in contradiction with the property of the inward normal half-line.
\proofend

\subsection{Graph estimates}

Let $\gamma$ be an embedded curve of class $C^{2,\alpha}$ in
$\R^2$, and let $p\in\Gamma:= Im(\gamma)$. Up to a rigid motion we
can assume that $p$ is the origin $O$ of $\mathbb{R}^2$ and
$\vec{\nu}(p)=(0,1)$. Let $\kappa$ be the curvature of $\Gamma$.
As $\Gamma$ is locally a graph, we have that around the origin
$\Gamma$ can be expressed as
\begin{equation}
\label{psi}
\psi(x) = (x,y(x)), \, \, \, \,
\end{equation}
with $y(0)=0$ and $y'(0)=0$ and then we have the following result.

\medskip

%%%xxx aclarar quien es "R" en el lema
\begin{lemma}\label{UGL}
If $|\kappa|\leq C$ in $\Gamma$ then, for any $p \in \Gamma$,
$\Gamma$ contains a graph (\ref{psi}) defined over the interval
$(-\varepsilon, \varepsilon)$ with $\psi(0)=p$. Here $\varepsilon$
depends only on $C$, and the functions $y(x)$, $y'(x)$ and $y''(x)$
are uniformly bounded in that interval.
\end{lemma}

{\it Proof.} As $|\kappa|\leq C$ and $|x|<\frac{1}{2C}$ by integrating the formula
\begin{equation}
\label{curvature}
\kappa(x) = \frac{d}{dx}\frac{y'}{\sqrt{1+(y')^2}} = \frac{y''}{(1+(y')^2)^{3/2}}
\end{equation}
we get
\[
\frac{|y'|}{\sqrt{1+(y')^2}} \leq C\, |x| < \frac{1}{2}.
\]
Observe now that the tangent vector $\vec{t}= \frac{1}{\sqrt{1+
(y')^2}}\,(1,y')$  satisfies $|\langle \vec{t}, (0,1) \rangle|<
1/2$. This inequality implies that the graph $\psi(x)$ can be extended to the interval $|x| < \varepsilon$ with $\varepsilon = 1/2C$. Moreover $|y'|$ is bounded in that interval in terms of
$\varepsilon$. We estimate the second derivative by using the
identity
\[
|{y''}| = |\kappa|\, (1+(y')^2)^{3/2} \,.
\]
and that proves the lemma. \proofend

\subsection{Harmonic overdetermined domains in the plane.}
When $f\equiv 0$, a classification of the domains of the plane
where problem (\ref{pr_bis}) is solvable is given in
\cite{traizet}. Assume that $\Omega$ is unbounded and $\partial
\Omega$ has a finite number of connected components; then, there
exist only three domains $\Omega$ where problem (\ref{pr_bis}) is
solvable (even for unbounded functions $u$!):
\begin{itemize}
\item the half-plane, \item the complement of a ball, and \item
the domain
\begin{equation}\label{HHP}
\Omega_{*} = \left\{(x,y) \in \R^2\,\,:\,\, |y| < \frac{\pi}{2} +
\cosh(x) \right\}\,
\end{equation}
that was first described in \cite{HHP}.
\end{itemize}

This correspondence gives in particular the following result.

\begin{lemma}\label{tra} (Corollary of Theorem 5 of \cite{traizet}). If $\Omega$ is a domain of the plane where
problem \eqref{pr_bis} can be solved for $f \equiv 0$, and the
boundary of $\Omega$ is unbounded and connected, then $\Omega$ is
a half-plane and $u$ is linear.
\end{lemma}

\section{Boundedness of the curvature}

The main result of this section is the following.

\begin{proposition}\label{curv}
Let $\Omega$ be an $f$-extremal domain with boundary unbounded and
connected, and $u$ a bounded solution to \eqref{pr_bis}. Then:

\begin{enumerate}

\item[i)] The curvature of $\partial \Omega$ is bounded.

\item[ii)] The $C^{2,\alpha}$ norm of the function $u$ is bounded
in $\overline{\Omega}$.

\end{enumerate}

\end{proposition}

{\it Proof.} If i) holds, Lemma \ref{regularity} implies a uniform
estimate of the $C^{2,\alpha}$ norm of $u$ near the boundary. In
the interior of $\Omega$, the $C^{2,\alpha}$ norm of $u$ is also
bounded due to interior regularity estimates (here we use in a
essential way the global boundedness of $u$). Therefore ii)
follows immediately.

\medskip

We now turn our attention to the proof of i). We recall that the
{\it accumulation set} of a sequence $F_n$ of subsets of
$\mathbb{R}^2$ is the closed set defined by
\[
Acc(F_n) = \{p \in \mathbb{R}^2: \ \exists\, p_n \in F_n\, \,
\textnormal{such that}\,\, p_n \to p\} \,.
\]

Let us suppose that $\partial \Omega$ has unbounded curvature, and
we will reach a contradiction. The proof uses a blow-up technique.

\medskip

{\bf Step 1: curvature rescaling.} Let $\kappa(q)$ denote the
curvature of $\partial \Omega$ at the point $q \in \partial
\Omega$. If $\kappa$ is unbounded, then there exists a sequence of
points $q_n \in \partial \Omega$ such that $|q_n|$ and
$|\kappa(q_n)|$ diverge to $+\infty$ increasingly. Let $I_n$ be
the connected component of $\partial \Omega \cap B_1(q_n)$
%%%xxx abierto mejor que cerrado. no afecta a lo demas y evita
%%%xxx problemas que podria traer el disco cerrado
containing $q_n$ and let $p_n =(x_n,y_n) \in I_n$ be the point
where the function
\[
p \to d(p, \partial B_1(q_n)) \, |\kappa(p)| = (1 - |p-q_n|)\,
|\kappa(p)|\, \, \, , \qquad p \in I_n
\]
attains its maximum, that clearly exists. We set
\[
r_n :=  d(p_n, \partial B_1(q_n))  = (1 - |p_n-q_n|)
\]
and
\[
R_n = r_n \,|\kappa(p_n)| \,.
\]

We have
\[
|\kappa(q_n)| \leq (1-|p_n-q_n|)\, |\kappa(p_n)| = r_n\,
|\kappa(p_n)| = R_n
\]
and then $R_n \to +\infty$. Since $r_n \leq 1$, we have also that
$|\kappa(p_n)|$ and $R_n/r_n$ diverge to $+\infty$. Consider the
transformation $T_n$ in $\mathbb{R}^2$ given by
\[
(x,y) \mapsto |\kappa(p_n)|\, (x-x_n, y-y_n) \,.
\]
Define $\Omega_n = T_n (\Omega)$. The image by $T_n$ of the balls
$B_{r_n}(p_n) \subset B_{1}(q_n)$ is given by the balls
$B_{R_n}(O)$, where $O$ is the origin of $\mathbb{R}^2$. If
$\kappa_n$ is the curvature of $\partial \Omega_n$, we have
clearly that
\[
\kappa_n = \frac{\kappa}{|\kappa(p_n)|}\,.
\]
Let $J_n = T_n(I_n)$. The function
\[
p \mapsto d(p, \partial B_{R_n/r_n}(O)) \, |\kappa_n(p)| =
(R_n/r_n - |p|)\, |\kappa_n(p)|\, \, \, , \qquad p \in J_n
\]
attains its maximum at $p=O$ and $|\kappa_n(O)| = 1$ for all $n$.
Let $R>0$. For $n$ large enough and $p \in J_n \cap B_R(O)$ we
have
\[
(R_n/r_n - R) \, |\kappa_n(p)| \leq  (R_n/r_n - |p|) \,
|\kappa_n(p)| \leq   (R_n/r_n - |O|) \, |\kappa_n(O)| = R_n/r_n\,.
\]
Then
\[
|\kappa_n(p)| \leq \frac{R_n/r_n}{R_n/r_n-R}
\]
for all $p \in J_n \cap B_R(O)$, and the curvature of
$\partial \Omega_n$ is uniformly bounded on compact sets.

\medskip

{\bf Step 2: existence of a limit curve.} Given $R>0$, define
$D_n(R)$ the connected component of $\Omega_n \cap B_R(O)$ which
has $O$ in its boundary, and $\Gamma_n(R)= \partial D_n(R) \cap B_R(O)$. By Lemma \ref{otro}, $\Gamma_n(R)$ is
connected.

Lemma \ref{UGL} implies the existence of $\delta>0$ such that
$\forall p \in \Gamma_n(R/2)$, the connected component of
$B_p(\delta) \cap \Gamma_n(R/2)$ passing through $p$ contains a
graph $Y_n$ of a function $y_n(x)$ defined on a segment of length
$\delta$. Moreover, Lemma \ref{regularity} implies that the
functions $y_n$ are of class $C^{2,\alpha}$ for all $\alpha \in
]0,1[$ and satisfy that their $C^{2,\alpha}$ norm is uniformly
bounded. Ascoli-Arzela's Theorem implies that a subsequence of
$y_n$ converges to a function $y_\infty \in C^{2,\alpha}(I_R(O))$
in the $C^{2,\alpha}$-topology, for all $\alpha \in ]0,1[$. A
prolongation argument allows to obtain
 $\Gamma_{n}$ and $D_n$ as the union of a subsequence of $\Gamma_n(R_n)$ and $D_n(R_n)$, where $R_n$ is chosen so that $R_n\to\infty$,  and a connected maximal
sheet $\Gamma_\infty$ of class $C^{2, \alpha}$ for all $\alpha \in
]0,1[$, such that $\Gamma_\infty$ belongs to the accumulation set
of $\{\Gamma_{n}\}$ and admits an arc-length parametrization
$\gamma_\infty(s)$ with $s \in \R$.

\medskip

{\bf Step 3: $\Gamma_\infty$ is proper.} If this was not the case,
there exists $p_n=\gamma_\infty(s_n) \in \Gamma_{\infty}$, where
$s_n$ is a divergent sequence and $p_n \to p\in \R^2$. By Lemma
\ref{UGL} and passing to the limit, there is a $\delta >0$ such
that each connected component of $\Gamma_{\infty} \cap
B_{\delta}(p)$ is a graph. Therefore we can choose $n$ so that
$\Gamma_{n} \cap {B_{\delta}(p)}$ has at least three connected
component which are curves passing through the consecutive points
$p_n,p_{n+1}$ and $p_{n+2}$. Enlarging $n$ if necessary, we can
assume that the distance of those components to $p$ is smaller
than $\delta/4$.

Now we can consider a connected component of $D_n \cap
B_{\delta}(p)$ with a boundary formed by two connected components,
both of them at a distance to $p$ smaller than $\delta/4$. Take $q
\in \partial\Omega_n$ with $|q-p| < \delta/4$. Then,
$B_{\delta/2}(q) \cap D_n$ gives a contradiction with Lemma
\ref{otro}.

\medskip

{\bf Step 4: $\Gamma_\infty$ is embedded.} By construction,
$\Gamma_\infty$ cannot have transversal self-intersections because
this would give rise to transversal self-intersections of
$\Gamma_n$ for $n$ large. But eventually $\Gamma_\infty$ could
have double tangential points, i.e. points $p$ such that there
exist $c_1 < c_2$ such that $p = \gamma_\infty(c_1) =
\gamma_\infty(c_2)$, and the tangent vectors to $\gamma_\infty$
satisfy $\vec{t}_\infty(c_1) = -\vec{t}_\infty(c_2)$. Let
$\gamma(s)$, $s \in \R$, be a parametrization of $\partial
\Omega$, $\vec{\nu}(s)$ be the unit normal vector of the curve
$\gamma(s)$ pointing to $\Omega$, and $\vec{\nu}_\infty$ its
induced limit unit normal on $\gamma_\infty$. Since
$\Gamma_\infty$ belongs to the accumulation set of $\{\Gamma_n\}$,
then the geometry of the curves $\gamma(s)$ and $\gamma_\infty(s)$
depend locally on the homotheties $T_n$ although the arc
parameters $s$ of these curves are not globally related.

We can suppose that the two values of $\vec{\nu}_\infty$ at $p$
are given by $(-1,0)$ and $(1,0)$. Moreover, as it is not possible
that all the points $\gamma_\infty(s)$ with $s \in [c_1, c_2]$ to
be double points, we can assume that $\gamma_\infty(s)$ is an
embedded curve in the open interval $c_1<s<c_2$ and there exist
\[
c_1 < c_3 < c_4 < c_2
\]
such that the angle between $\vec{\nu}_\infty(c_3)$ and
$\vec{\nu}_\infty(c_4)$, measured in the counterclockwise sense,
is strictly less than $\pi$.

\begin{figure}[!ht]
\centering {\scalebox{.6}{\input{ ./emb1.pstex_t}}} \caption{}
\label{fig9}
\end{figure}

By using the arc parameter $s$ of the curve $\gamma$ we get that
there exist four sequences $c_1^i$, $c_2^i$, $c_3^i$ and $c_4^i$,
$i=1,2,3,\cdots$, such that
\begin{itemize}
\item $c_1^1 < c_3^1 < c_4^1 < c_2^1 < c_1^2 < c_3^2 < c_4^2 <
c_2^2 < \cdots$ \item $|c_1^i-c_2^i|\to 0$, \item
$\vec{\nu}(c_1^i)$ converges to $(-1,0)$ and $\vec{\nu}(c_2^i)$
converges to $(1,0)$ as $i \to + \infty$, \item $\vec{\nu}(c_3^i)$
converges to $\vec{\nu}_\infty(c_3)$ and $\vec{\nu}(c_4^i)$
converges to $\vec{\nu}_\infty(c_4)$ as $i \to + \infty$.
\end{itemize}
In particular, for $i$ large enough we have that the angle between
$\vec{\nu}(c_3^i)$ and $\vec{\nu}(c_4^i)$, measured in the
counterclockwise sense, is strictly less than $\pi-\delta$ for
some $\delta>0$. This gives easily a contradiction with the
property of the normal inward half-line. In conclusion,
$\Gamma_\infty$ is embedded.

\medskip

{\bf Step 5: one of the connected components of
$\R^2\backslash \Gamma_\infty$ is contained in the accumulation set
$Acc(D_n)$.} The curve $\Gamma_{\infty}$ is properly embedded in
$\R^2$ and hence it separates $\R^2$ in two connected components.
Recall also that $D_n$ and $\Gamma_n$ are connected (the former by
definition, the latter by Lemma \ref{otro}). From these facts the
proof of Step 5 is elementary, and we denote by $\Omega_\infty$ the domain in $\R^2$ given by the connected component of
$\R^2\backslash\Gamma_\infty$ contained in the accumulation set
$Acc(D_n)$.

\medskip

{\bf Step 6: conclusion.} Clearly
$\Omega_{\infty}$ is a $C^{2,\alpha}$ domain and $\partial
\Omega_\infty = \Gamma_\infty$ with curvature $\kappa_\infty$
satisfying
\[
|\kappa_\infty(s)| \leq 1 = |\kappa_\infty(0)| \ \ {\rm for \ \ \
all} \ \ \ s \,.
\]
Define
\[
v_n(x,y) = |\kappa(p_n)|\, u\left( \frac{x}{ |\kappa(p_n)|} +x_n,
\frac{y}{|\kappa(p_n)|} + y_n \right),
\]

\[
f_n(x,y) = \frac{1}{ |\kappa(p_n)|}\, f\left(
\frac{1}{|\kappa(p_n)|}v_n (x,y) \right).
\]
Then $v_n$ solves the problem:
\begin{equation}
\left\{\begin{array} {ll}
\Delta v_n (x,y) +  f_n(x,y)  =0 & \mbox{in }\; \Omega_n\\
                        v_n > 0 & \mbox{in }\; \Omega_n\\
               v_n= 0 & \mbox{on }\; \pp \Omega_n \\
 \frac{\partial v_n}{\partial \vec{\nu}} = 1
&\mbox{on }\; \pp \Omega_n\,.
\end{array}\right.
\end{equation}
%%%xxx el antiguo p_0 podria estar a discancia cero del limite
Take $p_0 \in \R^2\setminus \Omega_\infty$, with distance to the
boundary greater or equal to a certain positive constant $r_0>0$
(this is possible by Step 5). Define $\psi$ in $B_R(p_0)$ as the
solution to the problem:
\begin{equation}
\left\{\begin{array} {ll}
\Delta \psi  = 4 \pi R^2 \delta_{p_0} - 1 & \mbox{in }\; B_R(p_0)\\
                        \psi= 0 & \mbox{on }\; \pp B_R(p_0)\,.
\end{array}\right.
\end{equation}
Here we denote by $\delta_{p_0}$ the Dirac delta measure centered
at $p_0$. Observe that the function $\psi$ is radially symmetric,
has a logarithmic singularity at $p_0$ and $\frac{\partial
\psi}{\partial \vec{\nu}}=0$ on $\partial B_R(p_0)$.

If $n$ is large enough, then $R_n > R$, $D_n \cap B_R(p_0)$ is
connected and the closure of $\Omega_\infty \cap B_R(p_0)$
coincides with $Acc(D_n \cap B_R(p_0))$.

We claim that $\|v_n\|_{C^{2,\alpha}}$ is bounded in $D_n \cap
B_R(p_0)$ for any $R$ fixed. For that we apply Green formula and
we obtain:
$$
\int_{D_n \cap B_R(p_0)} v_n - \psi\, f_n= \int_{D_n \cap
B_R(p_0)} (\psi\, \Delta v_n - \Delta \psi\, v_n)
$$
$$ =
\int_{\partial D_n \cap \overline{B_R(p_0)}} \left(\frac{\partial
\psi}{\partial \vec{\nu}} v_n - \psi \frac{\partial v_n}{\partial
\vec{\nu}}\right) = - \int_{D_n \cap \overline{B_R(p_0)}} \psi \,.
$$
Observe that the last term is uniformly bounded as it converges to
$$
\int_{\partial \Omega_\infty \cap \overline{B_R(p_0)}}\psi\,.
$$
Moreover $f_n$ is bounded in $L^{\infty}$. Hence we obtain that $
\int_{D_n \cap B_R(p_0)} v_n $ is bounded. Theorem 9.26 of
\cite{gilbarg} implies that $v_n$ is bounded in $L^{\infty}$ sense
on compact sets. Then we apply \cite{gilbarg}, Theorem 8.33 and
the comment that follows (see also Corollaries 8.35, 8.36 there),
to obtain that $v_n$ is bounded in $C^{1,\alpha}$ sense on compact
sets. In particular, $f_n$ is bounded in $C^{0,\alpha}$, always on
compact sets. Finally, Lemma \ref{regularity} yields the claim.

Then, by Ascoli-Arzela's Theorem $v_n$ converges in $C^{2,\alpha}$
sense (on compact sets) to a solution $v$ of the problem:
\begin{equation}
\left\{\begin{array} {ll}
\Delta v(x,y)  =0 & \mbox{in }\; \Omega_\infty\\
                        v > 0 & \mbox{in }\; \Omega_\infty\\
               v= 0 & \mbox{on }\; \pp \Omega_\infty \\
 \frac{\partial v}{\partial \vec{\nu}} = 1
&\mbox{on }\; \pp \Omega_\infty \,.
\end{array}\right.
\end{equation}
We apply now Lemma \ref{tra} and conclude that $\Omega_\infty$ is
a half-plane. But $\pp \Omega_\infty$ has a point (the origin $O$)
with curvature equal to $\pm 1$, and this yields the desired
contradiction. \proofend

\section{The case when $u$ is increasing in one variable}

The main result of this section is the following, that represents
the answer to our problem if $u$ is increasing in one variable.

\begin{proposition}\label{parc2}
Let $\Omega$ be a domain of $\mathbb{R}^2$ and suppose that $u =
u(x,y)$ is a solution of (\ref{pr_bis}) with $\frac{\partial
u}{\partial y}>0$ in $\overline{\Omega}$. Then $\Omega$ is a
half-plane and $u$ is parallel, that is, $u$
depends only on one variable.
\end{proposition}

\begin{remarks*} The same result is proved in \cite{Far-Val-0}, but under the
hypothesis that the domain is of class $C^3$. Our proof follows the one
in \cite{amb-cab, ber-caf-nir}.
\end{remarks*}

{\it Proof of Proposition \ref{parc2}.} Let $u_x$ and $u_y$ be the
derivatives of $u$ with respect to $x$ and $y$. Then, we can
define:
\[
\sigma = \frac{u_x}{u_y}, \ \ \ \ \ F= u_y^2\, \nabla \sigma.
\]
Then $\sigma$ is a function of class $C^{1,\alpha}$ in $\overline
\Omega$ and $F$ is just $C^{0,\alpha}$. We claim that
\begin{equation}\label{sigmadiv}
\nabla \cdot F=0 \ \mbox{ in } \Omega
\end{equation}
in the distributional sense. To see that, we multiply both sides
of $\Delta u + f(u)=0$ by a test function $\xi \in
C_0^{\infty}(\Omega)$ and integrate by parts, to obtain:
\begin{equation}
\label{sol1} \int_{\Omega} \left[ \langle \nabla \xi, \nabla
u\rangle - \xi f(u) \right] =0.
\end{equation}
Differentiating such equation with respect to $y$ we obtain
\[
\int_{\Omega} \left[ \langle \nabla {\xi_y}, \nabla u\rangle -
{\xi_y} f(u) \right] + \int_{\Omega} \left[ \langle \nabla \xi,
\nabla u_y\rangle - \xi f'(u)u_y \right] =0 \,.
\]
Then $v=u_y$ is a weak solution of
\begin{equation}
\label{sol2}
\Delta v + f'(u)v=0
\end{equation}
in $\Omega$ and the same holds for $v=u_x$.\\

Observe that $F= u_y \nabla u_x - u_x \nabla u_y$ and therefore,

$$ \int_{\Omega} F \cdot \nabla \xi = \int_{\Omega} u_y \nabla u_x \cdot \nabla \xi - u_x
\nabla u_y \cdot \nabla \xi= \int_{\Omega} u_y f'(u) u_x \xi - u_x
f'(u) u_y \xi=0$$
and (\ref{sigmadiv}) is proved.

\medskip

As $F$ is continuous in $\overline{\Omega}$ with $0$ divergence in
the distributional sense, by \cite{anze}[Remark 1.8] (see also
\cite{torres}[Theorem 7.2]), we have that the Divergence Gauss
theorem is valid in this framework. So, for any $\zeta \in
C_0^{\infty}(\mathbb{R}^2$ we have
\[
\int_\Omega \nabla \cdot (\zeta^2 \, \sigma\, F) = \int_{\partial
\Omega} \zeta^2\, \sigma\, \langle F, \vec{\nu}\rangle,
\]
where $\vec{\nu}$ is the inward normal vector about $\partial
\Omega$. Recall that in $\partial \Omega$, $\nabla u = \vec{\nu}$;
denoting by $\vec{e_1}, \vec{e_2}$ the vectors $(1,0)$ and
$(0,1)$, we have
\begin{eqnarray*}
\langle F, \vec{\nu}\rangle & = & \langle \nabla u_x, \vec{\nu} \rangle\, u_y
- \langle \nabla u_y\, , \vec{\nu} \rangle\, u_x  =
(\nabla^2 u) (\vec{\nu},\vec{e_1}) \, \langle \nabla u , \vec{e_2} \rangle -
(\nabla^2 u) (\vec{\nu},\vec{e_2}) \, \langle \nabla u , \vec{e_1} \rangle\\
 & = & (\nabla^2 u) (\vec{\nu},\vec{\nu}) \, \langle \vec{\nu} , \vec{e_1} \rangle\,  \langle \vec{\nu} , \vec{e_2} \rangle\,
 - (\nabla^2 u) (\vec{\nu},\vec{\nu}) \, \langle \vec{\nu} , \vec{e_2} \rangle\,  \langle \vec{\nu} , \vec{e_1} \rangle=0
\end{eqnarray*}
and then
\[
\int_\Omega \nabla \cdot (\zeta^2 \, \sigma\, u_y^2\, \nabla \sigma) =  0 \,.
\]
A simple computation gives
\[
\nabla \cdot (\zeta^2 \, \sigma\, F) = \zeta^2\, \sigma \, \nabla
\cdot F + 2\, \zeta\, \sigma \, u_y^2\, \langle \nabla \zeta,
\nabla \sigma \rangle + \zeta^2\, u_y^2\, |\nabla \sigma|^2
\]
and using \eqref{sigmadiv} we have
\[
\int_\Omega   \zeta^2\, u_y^2\, |\nabla \sigma|^2 = -2 \int_\Omega \zeta\, \sigma \, u_y^2\, \langle \nabla \zeta,  \nabla \sigma \rangle \,.
\]
From this last formula, using the H\"older inequality we obtain
\begin{eqnarray*}
\int_\Omega   \zeta^2\, u_y^2\, |\nabla \sigma|^2 &  \leq  & 2\, \left[\int_{\Omega} \zeta^2\, u_y^2\, |\nabla \sigma|^2 \right]^{1/2} \, \left[\int_{\Omega} u_y^2\, \sigma^2\,|\nabla \zeta|^2 \right]^{1/2} \,.
\end{eqnarray*}
By Proposition \ref{curv}, the gradient of $u$ is bounded, hence
so it is $u_x = u_y \sigma$. Therefore
\begin{equation}
\label{1}
\int_\Omega   \zeta^2\, u_y^2\, |\nabla \sigma|^2  \leq C_1\int_{\mathbb{R}^2}|\nabla \zeta|^2
\end{equation}
for some constant $C_1$. It is well known that in the plane there is a sequence of logarithmic cutoff functions
$\{\zeta_n\}_n\subset C_0^\infty(\mathbb{R}^2)$,  such that
\[
0 \leq \zeta_n \leq 1, \ \ \ \ \zeta_n = 1  \ \ in \ \ B_n(O) \ \
\ \lim_n \int_{\mathbb{R}^2}|\nabla \zeta_n|^2 = 0.
\]
Putting $\zeta = \zeta_n$ in (\ref{1}) and letting $n \to \infty$
we obtain
\[
\int_\Omega  u_y^2\, |\nabla \sigma|^2 = 0
\]
which means that $\sigma$ is constant, and then
\[
u_x (x,y) = C \, u_y(x,y)
\]
for a constant $C$. Then $\nabla u$ is normal to the vector
$(1,-C)$, and then $u$ is constant on every line parallel to that
vector,  i.e. $u$ is parallel. \proofend

\section{An $f$-extremal domain contains a tangent half-plane}

In this section we shall prove the following:

\begin{proposition}\label{tangent} Let $\Omega$ be an $f$-extremal domain in $\mathbb{R}^2$ whose boundary $\partial \Omega$ is unbounded and connected. Then $\Omega$ contains a half-plane $H$ such that $\partial \Omega$ and $\partial H$ are tangent.
\end{proposition}

Actually our result is stronger than stated in the above
proposition. In order to state and prove our results, we introduce
the concept of {\it limit direction} for the boundary of a domain.

\begin{definition} Let $\Omega$ be an unbounded domain in $\mathbb{R}^2$
with $\partial \Omega$ unbounded and connected, and let $P \in
\mathbb{R}^2$. We say that $v \in \mathbb{S}^1$ is a {\it limit
direction} for $\partial \Omega$ if there exists a sequence of
points $p_n \in \partial \Omega$ such that $|p_n| \to +\infty$ and
\[
\lim_{n \to +\infty} \frac{p_n-P}{|p_n-P|} = v \,.
\]
\end{definition}

Obviously, the set of limit directions is not empty and it does
not depend on the choice of the initial point $P$. We can fix the
coordinates of $\mathbb{R}^2$ in order that $O= (0,0) \in \partial
\Omega$, $\partial \Omega$ is tangent to the $x$-axis in $O$, and
the normal inward half-line $N(O)$ is the positive part of the
$y$-axis. Let $(\partial \Omega)_l$ and $(\partial \Omega)_r$ the
two components of $\partial \Omega \backslash \{O\}$, such that
$(\partial \Omega)_l$ near $O$ stays to the left of the $y$-axis,
and $(\partial \Omega)_r$ near $O$ stays to the right of the
$y$-axis.

\begin{definition}\label{lim_dir} We say that $v_l \in \mathbb{S}^1$ is a
{\it limit direction to the left} if there exists a sequence of
points $l_n \in (\partial \Omega)_l$ such that $|l_n| \to +\infty$
and
\begin{equation}\label{ln}
\lim_{n \to +\infty} \frac{l_n-O}{|l_n-O|} = v_l \,.
\end{equation}
We say that $v_r \in \mathbb{S}^1$ is a {\it limit direction to
the right} if there exists a sequence of points $r_n \in (\partial
\Omega)_r$ such that $|r_n| \to +\infty$ and
\begin{equation}\label{rn}
\lim_{n \to +\infty} \frac{r_n-O}{|r_n-O|} = v_r \,.
\end{equation}
\end{definition}

In particular, limit directions to the left or right are limit
directions. Moreover, there exist always at least one limit
direction to the left and one limit direction to the right. If
$v_l$ and $v_r$ are the limit directions respectively to the left
and to the right, let us denote by $\theta(v_l,v_r) \in [0,2\pi]$
the angle between $v_l$ and $v_r$ (measured from $v_l$ to $v_r$ in
the clockwise sense).

\begin{lemma}\label{ageqpi} Let $\Omega$ be un $f$-extremal domain
with boundary unbounded and connected. Let $v_l$ and $v_r$ be two limit
directions, respectively to the left and to the right. Then
\[
\theta(v_l,v_r)  \geq \pi \,.
\]
\end{lemma}

{\it Proof.} The proof uses an argument inspired by
\cite{esteban}. Let us suppose that $0 \leq \theta(v_l,v_r) <
\pi$, and let $l_n =(x_n^l, y_n^l)$ and $r_n = (x_n^r, y_n^r)$ be the
two sequences of points of $\partial \Omega$ (as in Definition
\ref{lim_dir}) for the limit directions $v_l = (x_{v_l}, y_{v_l})
\in \mathbb{S}^1$ and $v_r = (x_{v_r}, y_{v_r}) \in \mathbb{S}^1$.
After a suitable rotation and translation, we can suppose that $O
\in \partial \Omega$, $v_r \in (0, \pi/2 ]$, $v_l \in [\pi/2,
\pi)$. This means that, up to consider subsequences, we have that
$ y_n^l$, $y_n^r \to +\infty$.

\medskip Since $(\partial \Omega)_r$ is connected, it is possible to
replace the sequences $l_n$ and $r_n$ by other two sequences of
points such that $y_n^r = y_n^l$ for all $n \in \mathbb{N}$.
Hence, consider the segment $L_n$ that joins $l_n$ with $r_n$. It
is easy to see that the moving plane method applies to the part of
$\Omega$ that lies under $L_n$ with respect to horizontal lines.
Indeed, all connected componentes of $L_n^- \cap \Omega$ are
bounded, where $L_n^-=\{y \leq y_n^r\}$. Since this can be done for
all $n \in \mathbb{N}$, at the limit for $n \to +\infty$ we get
that the moving plane method applies to all $\Omega$ with respect
to horizontal lines. Then $\partial \Omega$ is a graph with
respect to the $x$-coordinates over an interval in $\mathbb{R}$,
and the solution $u$ of problem (\ref{pr_bis}) depends only on the
variable $y$ and is increasing in $y$. By Proposition \ref{parc2}
we conclude that $\Omega$ is a half-plane, but this is a
contradiction with the hypothesis that $0 \leq \theta(v_l,v_r) <
\pi$.

\proofend

Next lemma excludes from our study the case $\theta(v_l,v_r)
=\pi$.

\begin{lemma}\label{epi} Let $\Omega$ be $f$-extremal domain with boundary unbounded and connected, and $v_l$ and $v_r$ be two limit directions, respectively to the left and to the right, with $\theta(v_l,v_r)  = \pi$. Then
$\Omega$ is a half-plane and the bounded solution $u$ to problem
(\ref{pr_bis}) is parallel.
\end{lemma}

{\it Proof.} The proof of this lemma is inspired on the {\it
tilted moving plane} introduced in \cite{KKS} for constant mean
curvature surfaces. This procedure has also been applied to
elliptic problems in half-planes in \cite{sciunzi}, and to
overdetermined problems in \cite{Ros-Sic}.

Up to a suitable rotation $\mathcal{R}$ of $\mathbb{R}^2$ with
center in $O$ we can suppose that $v_l=(-1,0)$ and $v_r = (1,0)$.
Up to a translation, we can suppose that the origin $O$ of
$\mathbb{R}^2$ belongs to $\partial \Omega$, and $\partial \Omega$
intersect the $y$-axis transversally. Up to a reflection on the
$x$-axis, we can suppose that there exist $\delta>0$ such that
$\{(0,y)\, : \, \, 0\leq y\leq \delta\}$ stays in $\Omega$.

Consider $\Omega_1=\Omega\cap \{x>0\}$ and $\Omega_2=\Omega\cap
\{x<0\}$. Given a straight line $T$, for any $x\in\mathbb{R}^2$
and any subset $X\subset\mathbb{R}^2$ let $x'$ be the reflection
of $x$ about $T$ and $X'$ be the reflected image of $X$ about $T$.
Fix $\varepsilon>0$ small enough and consider the two families of
parallel straight lines
\[
T_a =\{ y= a\} \qquad \textnormal{and} \qquad T_{\varepsilon,a}=
\{y=-\varepsilon \, x+ a\}
\]
for $a \in \mathbb{R}$. Let $T=T_{\varepsilon,a}$ be an element of
the second family. For $a\geq0$ the line $T$ cut off from
$\Omega_1$ a bounded cap $\Sigma(T)$ defined as follows. As the
part of $\Omega_1$ below $T$ is made only by bounded connected
components (because $(1,0)$ is a limit direction of $\partial
\Omega$ to the right), it follows from the moving plane method
that the reflected image with respect to $T$ of the connected
components of $\Omega_1\cap\{y<-\varepsilon \, x +a\}$ is
contained in $\Omega$, except possibly for the component whose
boundary contains $O$. Let us denote this component by
$\Sigma(T)$. The portions of the boundary of $\Sigma(T)$ contained
in $T$, $x=0$ and $\partial\Omega$ will be denoted respectively by
$I$,
$J$ and $K$. Note that $O\in J\cap K$.\\
Let  $\Sigma'(T)$, $K'$, $J'$ and $O'$ be respectively the
symmetric image of $\Sigma(T)$, $K$, $J$ and $O$ with respect to
$T$. Define on the closure of $\Sigma'(T)$ the function $u'_T$
given by $u'_T (x) = u(x')$. At the beginning $\Sigma'(T)$ is
contained in $\Omega$ and $u'_T\leq u$ and we continue the process
while this
%while $\Sigma'(T)\subset \Omega$
occurs.

\begin{figure}[!ht]
\centering {\scalebox{.5}{\input{./tilted.pstex_t}}} \caption{}
\label{fig100}
\end{figure}

The process ends if we meet a first value
$a=a(\varepsilon)>0$ for which one of the following
events holds:
\begin{enumerate}
\item at an interior point, the reflected arc $K'$ touches the
boundary of $\Omega$, \item $K$ meets $T$ orthogonally, \item at a
point of $\Sigma'(T)\cup I$, the graph of the resulting function
$u'_T$ is tangent to the graph of the function $u$, \item $O'$
belongs to $\partial\Omega$, \item when restricted to the segment
$J'$, the graph of the resulting function  $u'_T$ is tangent at
some interior point to the graph of the function $u$.
\end{enumerate}
By the moving plane method, we deduce that each one of the first
three options implies that $K'\subset\partial \Omega$. Therefore
both events $(4)$ and $(5)$ are also true. We conclude that in
fact the process can be carried on either for all $a \geq 0$ or
until either event (4) or event (5) occurs for a first value
$a=a(\varepsilon)>0$. We can say that the process can be carried
on till $a$ reaches the limit value $a(\varepsilon)$, being
$a(\varepsilon) = +\infty$ if the process can be carried on for
all $a \geq 0$. Since $\partial \Omega$
intersect the $y$-axis transversally, we have that there exists a constant $C>0$ such that $a(\varepsilon)>C$.

\medskip Now take a sequence of $\varepsilon_i>0$ going to zero,
and repeat all the reasoning with $\varepsilon = \varepsilon_i$.
Let $a_1 \in [C,+\infty]$ be the limit of a subsequence of
$a(\varepsilon_i)$. If $a_1 = +\infty$, we conclude that the
moving plane method can be applied to $\Omega_1$ with respect to
all horizontal lines. If $a_1 \neq +\infty$, we conclude that the
moving plane method can be applied to $\Omega_1 \cap \{ y < a_1\}$
with respect to all horizontal lines and one of the two events
$(4)$ or $(5)$ for $T = T_{a_1}$. Moreover, since now $J$ is an
interval of $\Omega\cap\{x=0\}$, the value of $a_1$ depends only on
the behavior of $u$
restricted to $\Omega\cap\{x=0\}$.\\
Now repeat all the process for
$\Omega_2=\Omega\cap\{y < 0\}$ instead of
$\Omega_1$, with lines of positive slope defined
by $T^{*}_{\varepsilon,a}=\{y=\varepsilon x +a\}$.
We obtain that the moving plane can be applied either to $\Omega_2$ with respect to all horizontal lines (in this case we will define $a_2 = +\infty$), or to $\Omega_2 \cap \{ y < a_2\}$, for some $a_2 > 0$, with respect to all horizontal lines and one of the two events $(4)$ or $(5)$ holds for $T = T_{a_2}$. As it happens for $a_1$, the generalized number $a_2$ depends only on the behavior
of the solution $u$ along $\Omega \cap \{x=0\}$. From this last property, it follows that $a_1=a_2$. If $a_1 \neq +\infty$, the line $T_{a_1}$ satisfies that the reflected image of
$\Omega\cap \{y<a_1\}$ with respect to $T$ is contained
in $\Omega$, $u'_T\leq u$ and one of the assertions
$(1)$, $(2)$ or $(3)$ holds (at some point of the $y$-axis).
From the moving plane method we obtain that $\Omega$, and in particular $\partial \Omega$,
is symmetric with respect to $T$.

Now, we know that the origin $O \in \partial \Omega$ stays under
$T$. Since $\partial \Omega$ is connected, we have that $(\partial
\Omega)$ intersects $T$ and, since $\partial \Omega$ is symmetric
with respect to $T$, the vector $(1,0)$ would be a limit direction
to the left, which is not possible by Lemma \ref{ageqpi}. We
conclude that $a_1 = +\infty$, and then the moving plane method
can be applied to $\Omega$ with respect to all horizontal line.
Proposition \ref{parc2} concludes the proof.\proofend

From those two lemmas the proof of Proposition \ref{tangent} is
immediate.

\medskip

{\it Proof of Proposition \ref{tangent}}. By Lemmas \ref{ageqpi},
\ref{epi}, there are two possibilities: either $\Omega$ is a
half-plane or $\theta(v_l,v_r)
> \pi$ for any limit directions to the left and right. In this last case, assume that $O \in \partial
\Omega$, and make a convenient rotation so that:
$$v_r<0 \mbox{ and } v_l > \pi,$$
for any limit directions to the right and left, respectively. Then
$\partial \Omega \cap \{ y \geq -1\}$ is non-empty and compact.
Therefore there exists $c\geq 0$ so that $H=\{y \geq c\}$ is the
claimed half-plane. \proofend

\section{Building a parallel solution starting from an $f$-extremal domain}

The main result of this section is the following:

\begin{proposition}\label{HHH}
There exists a sequence of points $q_n \in
\partial \Omega$ satisfying that:

\begin{enumerate}

\item $|q_n| \to + \infty$ and $\frac{q_n}{|q_n|}\to v \in
\mathbb{S}^1$ for some direction to the right $v$.

\item If $T_n$ is the translation in $\R^2$ that moves $q_n$ to
the origin, then $\Omega_n = T_n(\Omega)$ converges to the
half-plane $$\Omega_\infty=\{ p \in \R^2:  \langle v^{\perp} , p\rangle
>0\}.$$ Here $v^{\perp}$ denotes the vector obtained by rotating $v$ of angle $\pi/2$ measured in the counterclockwise sense.
Moreover, the sequence of functions $u_n(x,y) = u((x,y)-q_n)$
converges to a bounded parallel solution of (\ref{pr_bis}) in
$\Omega_\infty$.
\end{enumerate}

\end{proposition}

\begin{remark} \label{puesclaro} An analogous statement is true
for a certain direction to the left $\tilde{v}$. \end{remark}

\medskip

{\it Proof.} As always, we can assume that the origin $O$ belongs
to $\partial \Omega$.  Observe that the set of the limit
directions to the right is closed. Moreover, it is not the whole
$\mathbb{S}^1$ because $\Omega$ contains a half-plane. Then, we
can choose $v = e^{i\theta}$ a limit direction to the right such
that $e^{i(\theta-\epsilon)}$ is not a limit direction to the
right for any $\epsilon \in (0, \epsilon_0)$. Up to a rigid
motion, we can assume that $v=(1,0)$.

\medskip

Take $\epsilon$ small enough. Consider the sector of $\R^2$ given by
\[
C_\epsilon = \{(x,y)\in \R^2:\ |y| \leq \epsilon x\, \} \,.
\]

By the choice of $v$ we know that $(\partial \Omega)_r \cap
C_{\epsilon}$ is unbounded but the part of $(\partial \Omega)_r$
that lies under $C_\epsilon$ is compact. If $p=(x^p,y^p) \in
(\partial \Omega)_r$ we define also the sector of $\R^2$ given by
\[
G_{p,\epsilon} = \left\{(x,y):\  y \leq y^p - 2 \epsilon |x-x^p| \right\}.
\]

Choose $p_\epsilon = (x_\epsilon,y_\epsilon)$ such that
\begin{itemize}
\item $p_\epsilon \in C_\epsilon \cap (\partial \Omega)_r$, \item
the distance of $p_\epsilon$ to the origin is bigger than
$1/\epsilon^2$
\end{itemize}

 Observe that $(\partial \Omega)_r \cap G_{p_\epsilon, \epsilon}$ is compact
 and contained in $C_{\epsilon}$ if $\epsilon$ is sufficiently small.
 In particular there exists a point $q_\epsilon \in (\partial \Omega)_r \cap G_{p_\epsilon} \cap C_\epsilon$
 minimizing the function $(x,y) \mapsto y$ (see Figure \ref{fig10}).
 Such value $q_\epsilon$ satisfies that:

 \begin{itemize}
\item $|q_\epsilon| \to +\infty$ as $\epsilon \to 0$. \item
$(\partial \Omega)_r \cap G_{q_\epsilon,\epsilon}=
\{q_\epsilon\}.$
\end{itemize}

Now let $D_\epsilon$ be the connected component of
$B_{1/\sqrt{\epsilon}}(q_\epsilon) \cap \Omega$ containing
$q_\epsilon$ in its boundary. Observe that $D_\epsilon$ is above
the sector $G_{q_\epsilon}$. We do a translation $T_\epsilon$ in
$\R^2$, moving $q_\epsilon$ to the origin $O$, and we set
$D'_\epsilon = T_\epsilon(D_\epsilon)$.

\begin{figure}[!ht]
\centering {\scalebox{.5}{\input{./ekeland.pstex_t}}} \caption{}
\label{fig10}
\end{figure}

\medskip

We now make $\epsilon$ converge to $0$. By Proposition \ref{curv},
the curvature of $\Omega$ and the $C^{2,\alpha}$ norm of $u$ in
$\overline{\Omega}$ are bounded. Following the arguments in the
proof of Proposition \ref{curv} (see in particular Steps 2, 3, 4,
5) we have that the domains $D'_\epsilon$ converges to an
$f$-extremal domain with boundary unbounded and connected
$\Omega_\infty$. Since $G_\epsilon = T_\epsilon(G_{q_\epsilon})$
converges to the half-plane $\{y>0\}$, the domain $\Omega_\infty$
is contained in a half-plane. By Lemma \ref{epi} $\Omega_{\infty}=
\{y >0\}$, and then the sequence $u_n$ converges to a bounded
parallel solution. \proofend

\section{Existence of solutions in balls and asymptotic properties}

The main result of this section is the following:

\begin{proposition} \label{teo-ball}
Assume that there exists a solution of the problem:
\begin{equation} \label{eq-ODE} \left \{ \begin{array}{l}
\varphi''(y) + f(\varphi(y))=0 \\ \varphi(0)=0, \ \varphi'(0)= 1, \\
\displaystyle \lim_{t \to +\infty} \varphi(y)= L>0.\end{array} \right.
\end{equation}
Then, there exists $R_0>0$ such that for any $R>R_0$ the problem:
\begin{equation}\label{eq-ball} \left \{ \begin{array}{lr} \Delta u + f(u)=0 & x \in B_R(O),
\\ u>0, & x \in B_R(O), \\ u=0, & x \in \partial B_R(O).\end{array} \right. \end{equation}
admits a radially symmetric solution $u_R$. Moreover, $u_R$ has
the following asymptotic behavior:

\begin{enumerate}

\item[i)] $u_R < L$ and for any $\rho \in (0,1)$, $u_R|_{B_{\rho
R}(O)}$ converges uniformly to $L$ as $R \to +\infty$.

\item[ii)] The functions $v_R(z)= u_R(z-(0,R))$ converges to
$u(x,y)=\varphi(y)$ locally in compact sets of $H=\{ y > 0\}$.

\end{enumerate}

\end{proposition}

\medskip

\begin{remark} Actually, it will be clear from the proof that if $f(0)\geq 0$ there exist solutions
of \eqref{eq-ball} for any $R>0$. Instead, if $f(0)<0$ such
existence result is limited to large values of the radius.
\end{remark}

\medskip

In order to prove Proposition \ref{teo-ball}, we need some
preliminary work. First, we show that the existence of the ODE
\eqref{eq-ODE} is equivalent to some properties on $f$ and its
primitive, denoted by:
\begin{equation} \label{primitive} F(u)=\int_0^u f(s) \, ds\,.
\end{equation}

\medskip

\begin{lemma}\label{lem aux} The following two assertions are equivalent:

\begin{enumerate}

\item[i)] There exists a solution to \eqref{eq-ODE}.

\item[ii)] $f(L)=0$ and $F(L)=1/2 > F(u)$ for all $u \in [0,L)$.

\end{enumerate}
Moreover, in such case, there exists a sequence $\mu_n<L$, $\mu_n
\to L$ such that $F(\mu_n)> F(u)$ for all $u \in [0, \mu_n)$.

\end{lemma}

{\it Proof: i) $\Rightarrow$ ii).}  The limit at infinity of
$\varphi$ in \eqref{eq-ODE} implies that $f(L)=0$. Moreover, let
us recall that the Hamiltonian:
$$H = \frac{1}{2} (\varphi')^2 + F(\varphi),$$
is a constant in $y$. Observe that $\varphi'(y) \to 0$ if $y \to +
\infty$, so that such constant is nothing but $F(L)$. Moreover,
replacing $y=0$ we obtain the exact value of $F(L)$:
\begin{equation} \label{hola} F(L)= \frac{1}{2} \varphi'(0)^2 +
F(\varphi(0))=\frac{1}{2}. \end{equation}
Moreover, it is easy to observe that $\varphi'(y)>0$ for any $y \geq 0$. Then,
$$F(L)=\frac{1}{2} \varphi'(y)^2 + F(\varphi(y)) > F(\varphi(y))
\ \forall \ y \in [0,+\infty).$$

\bigskip {\it ii) $\Rightarrow$ i).} In the phase space, let us
consider the level set:
$$ \mathcal{C}= \{(\varphi, \varphi') \in [0,+\infty)^2:\ \frac{1}{2} (\varphi')^2 +
F(\varphi)=1/2\}.$$
This is a smooth curve for any $\varphi'>0$ as the Implicit
function theorem shows. Moreover, for any $\varphi \in [0,L]$,
there exists a unique $\varphi' \geq 0$ such that $(\varphi,
\varphi') \in \mathcal{C}$. Observe that $\varphi'=1$ if
$\varphi=0$ and $\varphi'=0$ if and only if $\varphi=L$.
Then, the solution of the Initial Value Problem:
\begin{equation} \left \{ \begin{array}{l}
\varphi''(y) + f(\varphi(y))=0 \\ \varphi(0)=0, \ \varphi'(0)=
1,\end{array} \right.
\end{equation}
has image in $\mathcal{C}$. Since $\varphi'>0$ for any $\varphi
\in (0,L)$, the image of the solution contains all $\mathcal{C}$
except, eventually, the point $(L, 0)$.

We now show that $\lim_{t \to +\infty} \varphi(t)=L$. Otherwise,
$\varphi$ arrives to the value $L$ at a certain time $t$, and
$\varphi'(t)=0$. However, since $f(L)=0$, $L$ is an equilibrium of
the ODE, and this gives a contradiction with the uniqueness of the
solution for the initial value problem.

\medskip

Observe that the last assertion of Lemma \ref{lem aux} would be
immediate if $f$ where positive below the value $L$, and actually
we would have a continuum of values satisfying such condition. In
general, though, $f$ could change infinitely many times below $L$.
Define
$$m_n= \max\left\{ F(x):\ x \in \left[0,L- \frac 1 n\right]\right\}, \ \mbox{and } \mu_n
= \min\left\{x \in \left[0,L- \frac 1 n\right]: \ F(x)=m_n\right\}.$$
By the definition of $\mu_n$, $F(\mu_n)=m_n > F(x)$ for all $x \in
[0,\mu_n)$.
We now show that $\mu_n \to L$. Otherwise, we could pass to a
subsequence (still denoted by $\mu_n$) such that $\mu_n \to \mu
<L$. Then, $F(\mu) \leftarrow F(\mu_n) = m_n \to F(L)$, which
implies that $F(\mu)=F(L)$, contradicting ii). \proofend

Our intention is now to settle the problem variationally. For
that, we need to truncate the function $f$ conveniently for $u<0$
and $u>L$.
Given $\delta>0$, we define:
$$ \tilde{f}(u) = \left \{ \begin{array}{lr} 0 & \mbox{ if } u \geq L , \\
f(u) & \mbox{ if } u \in [0,L], \\ f(0)(1 + \frac{u}{\delta}) & \mbox{ if }
u \in [-\delta, 0], \\ 0 & \mbox{ if } u \leq -\delta. \end{array} \right.$$
Accordingly, we define $\tilde{F}(u) = \int_0^u \tilde{f}(s) \,
ds$. Observe that for $u \leq -\delta$, $\tilde{F}(u)= - f(0)
\delta/2$.

We now fix $\delta>0$ so that
\begin{equation} \label{delta} F(u) > |f(0)|
\delta/2 \  \ \forall u \in [L-2\delta, L]. \end{equation}
It is then clear that:
\begin{equation} \label{cond-F} 1/2=\tilde{F}(L) > \tilde{F}(u) \ \forall u <L, \
\mbox{ and } \tilde{F}(\mu_n) > \tilde{F}(u)\ \forall u < \mu_n
\end{equation}
where $\mu_n$ is given by Lemma \ref{lem aux}, c), and we consider
only the terms of the sequence so that $|\mu_n - L|< \delta$.
With this truncation, our aim is to find solutions of the problem:
%%%xxx cambiar la frase. ahora suponemos u>0 y en el lema demostramos u(z) esta entre -\delta y L
\begin{equation}\label{eq-trunc} \left \{ \begin{array}{lr} \Delta u + \tilde{f}(u)=0 & \mbox{ in } B_R(O),
\\ u>0, & \mbox{ in } B_R(O), \\ u=0, & \mbox{ in } \partial B_R(O).\end{array} \right. \end{equation}

\medskip

\begin{lemma}\label{mp} Let $u$ be a solution of
\eqref{eq-trunc}. Then
\begin{equation} \label{u-trunc} \left \{ \begin{array}{lr} u(z) \in (0, L) & \mbox{ if } f(0)\geq  0, \\ u(z) \in
(-\delta, L) & \mbox{ if } f(0) <   0, \end{array} \right. \
\forall \ z \in B_R(O).\end{equation}
\end{lemma}

{\it Proof.} First let us show that $u(z)\leq L$ for any $z \in
B_R(O)$. Otherwise, assume that $\max u= u(z_0)>L$. Let
$\Omega=\{z \in B_R(O):\ u(z) >L\}$. Clearly $u$ is harmonic in
$\Omega$ and attains a maximum in its interior, which is
impossible. In the same way we can prove that $u(z)\geq 0$ (if
$f(0)\geq 0$) or $u(z) \geq - \delta$ (if $f(0)<0$).

We now show the strict inequality. Otherwise, assume that $\max u
= L$. Observe also that the constant function $L$ is a solution of
$\Delta u + f(u)=0$. Therefore both solutions are in contact, and
this is in contradiction with the maximum principle. \proofend

Let us define the energy functional:
$$I_R:H_0^1(B_R(O)) \to \R, \ I_R(u)= \int_{B_R(O)} \frac 1 2
|\nabla u|^2 - \tilde{F}(u).$$ Here $H_0^1(B_R(O))$ denotes the
closure of the space $C_0^{\infty}(B_R(O))$ with the usual Sobolev
norm
$$\|u \| = \left ( \int_{B_R(O)} |\nabla u|^2 + u^2 \right)^{1/2}. $$
Clearly \eqref{eq-trunc} is the Euler-Lagrange
equation of the functional $I_R$. The following lemma establishes
the existence of a minimum for $I_R$ and, therefore, a solution
for \eqref{eq-trunc}.

\begin{lemma}\label{variational} For any fixed $R>0$, the
functional $I_R$ attains its minimum at a radially symmetric
function $u_R$. \end{lemma}

{\it Proof.} This is quite standard. Observe that since
$\tilde{F}$ is continuous and bounded, the energy functional $I_R$
is coercive and weakly lower semi-continuous. From this we obtain
the existence of a minimizer $u_R$. By making use of the Schwartz
rearrangement (see for instance \cite{kawohl}), we can assume that
$u_R$ is radially symmetric.
\proofend

Observe that if $f(0) \geq 0$, by Lemma \ref{mp} we already have a
solution of our problem \eqref{eq-ball}. Instead, if $f(0)<0$ we
still need to show that $u_R$ is positive. But, before, let us
give some energy estimates on $u_R$.

In what follows, we denote by $A(p;R_1,R_2)$ the annulus of radii
$R_1<R_2$.

\medskip

\begin{lemma} \label{lem-energy} There exists $C>0$
independent of $R$ so that:
\begin{equation} \label{estim-energy} -\frac 1 2 \pi R^2 \leq I_R(u_R)
\leq -\frac 1 2 \pi R^2 + C R, \end{equation}
\begin{equation} \label{estim-F} \frac 1 2 \pi R^2 \geq \int_{B_R(O)}
\tF(u_R) \geq \frac 1 2 \pi R^2 - C R. \end{equation}
\end{lemma}

{\it Proof.} Taking into account \eqref{cond-F}, we have
$$I_R(u_R) =
\int_{B_R(O)} \frac 1 2 |\nabla u_R|^2 - \tilde{F}(u_R) \geq -
\int_{B_R(O)} \tilde{F}(u_R) \geq - \frac 1 2 \pi R^2.$$
From this we obtain the first inequality of \eqref{estim-energy}
and \eqref{estim-F}.

For the second inequality, let us define $\phi_R \in
H_0^1(B_R(O)$,
$$ \phi_R(|z|) = \left \{ \begin{array}{lr} L & |z| \leq R-1, \\
L(R-|z|) & |z| \in [R-1, R]. \end{array} \right.$$
We now estimate $I_R(\phi_R)$. The gradient term can be estimated
as:
$$ \int_{B_R(O)} |\nabla \phi_R|^2 = 2 \pi \int_{R-1}^R
\phi_R'(r)^2 r \, dr \leq C R.$$
In order to estimate the term $\int_{B_R(O)} \tilde{F}(\phi_R)$,
we split it into two terms:
$$\int_{B_{R-1}(O)} \tilde{F}(\phi_R)= \frac 1 2 \pi (R-1)^2 \geq \frac 1 2 \pi
R^2- C R,$$
$$\int_{A(0;R-1,R)} \tilde{F}(\phi_R) \geq - C R.$$
In the last estimate we have just used the boundedness of
$\tilde{F}$.
The above estimates imply that $I_R(\phi_R) \leq -\frac 1 2 \pi
R^2 + C R$. Since $I_R(u_R) \leq I_R(\phi_R)$, we conclude
\eqref{estim-energy}.
Finally,
$$ - \int_{B_R(O)} \tilde{F}(u_R) \leq I_R(u_R) \leq -\frac 1 2 \pi R^2 +
C R,$$
and \eqref{estim-F} follows.
\proofend

Next lemma is devoted to show the asymptotic behavior of
$u_R$.

\medskip

\begin{lemma} \label{prop} The following assertions hold:

\begin{enumerate}

\item[a)] For any fixed $\rho <L$, there exists $C=C_\rho$
independent of $R$ so that:
$$\Omega_\rho = \{ z \in B_R(O): u_R(z) < \rho\} \subset A(0;R-C_\rho,
R).$$

\item[b)] There exists $R_0>0$ such that $u_R$ is positive for $R
\geq R_0$.

%\item[c)] $\lim_{R\to + \infty} \frac{\partial u_R}{\partial \nu}
%= - \alpha$, where $K$ is as in \eqref{eq-ODE}.
\end{enumerate}
\end{lemma}

{\it Proof.} The proof of a) will be made in two steps.

\medskip

{\bf Step 1.} For any fixed $\rho <L$, there exists $C=C_\rho$
independent of $R$ so that $|\Omega_\rho| \leq C_\rho R$. Indeed,
$$ \int_{B_R(O) \setminus \Omega_\rho} \tF(u_R) \leq \frac 1 2 (\pi R^2
- |\Omega_\rho|),$$
$$ \int_{\Omega_\rho} \tF(u_R) \leq \max\{\tF(x): x < \rho\}
|\Omega_\rho| = \left(\frac 1 2-\varepsilon\right) |\Omega_\rho|,$$ where
$\varepsilon= \frac 1 2-\max\{\tF(x): x < \rho\}>0$ by
\eqref{cond-F}.
Adding both terms, we get:
$$ \int_{B_R(O)} \tilde{F}(u_R) \leq \frac 1 2 \pi R^2 - \varepsilon
|\Omega_\rho|,$$ and Step 1 follows from \eqref{estim-F}.

\medskip

{\bf Step 2.} Let us fix $R>0$ and $\mu=\mu_n$ one of the elements
of the sequence in Lemma \ref{lem aux} satisfying \eqref{cond-F}.
Then $\Omega_{\mu} = \{ z \in B_R(O): u_R(z) < \mu\}$ is
connected.

Observe that $\Omega_{\mu}$ always has a connected
component touching the boundary $\partial B_R(O)$. Suppose by
contradiction that it has an interior connected component too,
denoted by $U$. Then, $u_R(z) < \mu$ for $z \in U$ and $u_R(z)=
\mu$ if $z \in
\partial U$.

Define:
$$ v(z)= \left \{ \begin{array}{lr} u_R(z) & z \notin U, \\ \mu & z
\in U. \end{array} \right. $$
\begin{figure}[!ht]
\centering {\scalebox{.5}{\input{./david.pstex_t}}} \caption{}
\label{fig3}
\end{figure}

Clearly, $v \in H_0^1(B_R(O))$ and $\int_{U} |\nabla u_R|^2 \geq
\int_{U} |\nabla v|^2=0$. Moreover, taking into account
\eqref{cond-F},
$$\int_U \tF(u_R) \leq \int_U \tF(\mu) = \int_U \tF(\mu).$$
Therefore $I_R(v) < I_R(u_R)$, a contradiction that proves Step 2.

Step 1 and 2 readily imply a). Indeed, given $\rho<L$, take
$\mu=\mu_n \in (\rho, L)$ one of the elements of the sequence. Since $\Omega_\mu$ satisfies
the statements of Step 1 and 2, $\Omega_\mu \subset A(0, R-C, R)$
for some positive constant $C$. But
$\Omega_\rho \subset \Omega_{\mu}$, concluding the proof.

\bigskip

We now turn our attention to assertion b). The case $f(0) \geq 0$
is clear from Lemma \ref{mp}, so let us consider the case $f(0)
<0$. Suppose that there exists $r_0 \in [0,R)$ with
$u_R(r_0)=-\delta_R \leq 0$, $u_R'(r_0)=0$. By Lemma \ref{mp},
$\delta_R \in [0, \delta)$. Moreover, by a) we have that $r_0 \in
(R-C, R)$ for some positive $C>0$ independent of $R$.

Define $v(z)= u_R(z) + \delta_R$, which is a solution of the
problem:
$$ \left \{ \begin{array}{lr} \Delta v + g(v)=0 & \mbox{ in
}B_{r_0}(O), \\ v=0 & \mbox{ in } \partial B_{r_0}(O), \\
\frac{\partial v}{\partial \nu}=0 & \mbox{ in } \partial
B_{r_0}(O)
\end{array} \right. $$
where $g(t)= \tf(t- \delta_R)$. We now apply the Pohozaev identity
(see \cite{Struwe}[Chapter III, Lemma 1.4]) to the previous
problem, to obtain that
\begin{equation} \label{cip} \int_{B_{r_0}(O)} G(v)=0, \end{equation}
with $G(t)= \int_0^t g(s) \, ds = \tF(v-\delta_R) -
\tF(-\delta_R)$.

We will show now that this is impossible if $R$ is sufficiently large.
Indeed, take $\Omega_\mu$ the set defined in Step 2. Then,
$$ \int_{B_{r_0}(O)\setminus \Omega_{\mu}} G(v) = \int_{B_{r_0}(O)\setminus
\Omega_{\mu}} \tF(u_R- \delta_R) - \tF(-\delta_R).$$
Now, $\tF(-\delta_R) \leq |\tF(-\delta)| = |f(0)|
\frac{\delta}{2}$. Moreover, in $\Omega_\mu$, $u_R- \delta_R \geq
\mu - \delta_R \geq L - 2 \delta$. By \eqref{delta}, we conclude
that $\tF(u_R- \delta_R) - \tF(-\delta_R)>c>0$ for any $z \in
\Omega_\mu$. Then,
$$ \int_{B_{r_0}(O)\setminus \Omega_{\mu}} G(v) \geq c | B_{r_0}(O) \setminus \Omega_{\mu}| \geq c' R^2.$$
Moreover,
$$ \left |\int_{\Omega_\mu} G(v) \right | \leq \int_{A(0;R-C,R)}
|G(v)| = O(R),$$ and hence \eqref{cip} cannot hold for large $R$.
%In order to prove c), we use again Pohozaev identity
%(\cite{Struwe}[Chapter III, Lemma 1.4]):
%
%$$ 2 \int_{B_R(O)} F(u_R) = \frac{1}{2} R \int_{\partial B_R(O)}
%\left ( \frac{\partial u_R}{\partial \nu}\right )^2$$
%
%Observe that
%
%$$\frac{1}{2} R \int_{\partial B_R(O)} \left ( \frac{\partial
%u_R}{\partial \nu}\right )^2 = \pi R^2 \left ( \frac{\partial
%u_R}{\partial \nu}\right )^2.$$
%
%Then, by Lemma \ref{lem-energy},
%
%$$\lim_{R \to +\infty} \frac{\partial
%u_R}{\partial \nu} = \pm \sqrt{2 \frac 1 2}.$$ Since $u_R$ is positive,
%then $\frac{\partial u_R}{\partial \nu}$ must be negative, and we
%conclude by \eqref{hola}.
\proofend

\medskip

We are now able to prove Proposition \ref{teo-ball}.

\medskip

{\it Proof of Proposition \ref{teo-ball}} With the previous results,
we just need to prove iii). Take $R_n \to +\infty$, $n \in
\mathbb{N}$, $v_n=u_n(z-(0,R_n))$. We first show that in any
compact set of $H=\{(x,y) \in \R^2:\ y
>0\}$, $v_n$ is bounded in with respect to the $C^{2,\alpha}$ norm. We use a bootstrap argument in two steps: $f(u_n)$ is a uniformly
bounded function, and then $u$ is of class $C^{1,\alpha}$. Then,
$f(u_n)$ is a Lipschitz function, and we repeat the argument with
$C^{2,\alpha}$ regularity.

As a consequence, $v_n$ converges (up to a subsequence) to a
solution of the problem $\Delta v + f(v)=0$ defined in $H$. This
convergence is $C^{2,\alpha}$ in compact sets of $H$, with
$0<\alpha<1$. We now claim that $v$ is parallel.

Take $p=(x, y) \in H$. We denote by $\rho_{n}$ its distance to the
center of the ball $(0,R_n)$, that is, $\rho_{n}= \sqrt{x^2 +
(R_n-y)^2}$. Since $u_n$ is radially symmetric, then
$v_n(p)=v_n(0,R_n-\rho_{n})$. Observe now that that $R_n-\rho_{n}
\to y$. Therefore $v_n(p) \to v(0,y)$, which is independent of
$x$.

We now prove that $v(0)=0$. With the previous information, we can
consider the convergence of sequence $v_n(r)= u_n(r-n)$, which
solves:
$$ v_n''(r) + \frac{v_n'(r)}{r} + f(v_n(r))=0,\ \ v_n(0)=0.$$
If we consider that equation in $r \in [0,1]$, it is easy to show that it converges in $C^{2,\alpha}$ sense to $v(r)$.
In particular, $v(0)=0$.

Finally, we will show that $v=\varphi$ given in \eqref{eq-ODE} by
showing that $\lim_{r\to + \infty} v(r)=L$. Observe that Lemma
\ref{mp} implies that actually $v(r) \leq L$ for any $r \in
(0,+\infty)$.
Fix now $\rho>0$ and take $C>0$ as given by Lemma
\ref{prop}, a). Then, for any $r \in (C, 2 R_n-C)$ we have that
$v_n(r) \geq \rho$. As a consequence, $v(r) \geq \rho $ for any $r
\in (C, +\infty)$, which implies that $\lim_{r\to + \infty}
v(r)=L$.
Therefore we have proved the convergence of an adequate
subsequence. The uniqueness of the limit implies that actually the
whole sequence converges.\proofend

\section{Proof of the main theorem}

In order to conclude the proof of our main theorem, we recollect
the information from the previous sections.

\medskip

From Proposition \ref{tangent} we know that $\Omega$ contains a
half-plane $H$ internally tangent to $\partial \Omega$. Moreover,
from Lemmas 5.4, 5.5 we can assume that $\theta (v_l,v_r)> \pi$
for any directions to the left and right $v_l, v_r$. We can
suppose that $H$ is the half-plane $\{y>0\}$ and the interior
tangent point with $\partial \Omega$ is the origin. Moreover, by
Proposition \ref{HHH} there exists an unbounded sequence $q_n \in
\partial \Omega$ such that,
doing translations in $\R^2$ that move $q_n$ to the origin we get
a sequence of domain $\Omega_n$ converging to a limit half-plane
$\Omega_{\infty}$, and a sequence of functions $u_n$ converging to
a parallel solution. Recall also that $\frac{q_n}{|q_n|} \to v $
and $\Omega_{\infty} = \{p \in \R^2: \langle p , v^{\perp}\rangle
>0\}$. By making a small rotation and reflection, if necessary, we can assume
that $v= e^{i\theta}$, $\theta \in (-\pi/2,0)$.

In section 7 we proved that for every $R$ large enough there
exists a radial solution $u_R$ to the problem \eqref{eq-ball},
such that as $R \to +\infty$

\begin{enumerate}

\item[i)] $u_R < L$ and $u_R|_{B_{\rho R}(O)}$ converges uniformly
to $L= \lim_{y\to+\infty} \varphi$ for any $\rho \in (0,1)$,

\item[ii)] the functions $v_R(z)= u_R(z-(0,R))$ converges to
$\varphi$ locally in compact sets of $H$, where $\varphi$ is a
solution of \eqref{eq-ODE}.

\end{enumerate}

We are now able to prove our main result.

\medskip

{\it Proof of Theorem \ref{Tmain}.} Let $R$ large enough and
consider the solution $u_R$. Since the parallel solution
$u_\infty$ is obtained as limit of a sequence $u_n$ of
translations of the function $u$ in $\Omega$, we get that there
exists a point $p \in \Omega$ such that the ball $B_R(p)$ is
contained in $\Omega$ and the graph of the function $u_R$ defined
in $B_R(p)$ stays under the graph of the function $u$. Moreover,
$p$ can be chosen so that $|p-q_n|< 2R$ with $q_n$ an element of
the sequence described above.

Now, we claim that we can move the ball $B_R(p)$ inside $\Omega$
till it reaches the position of the ball $B_R(q)$ with $q=(0,R)$.
Observe that the graph of the function $u_R$, during the motion,
cannot touch the graph of the function $u$ by the maximum
principle.

\begin{figure}[!ht]
\centering {\scalebox{.5}{\input{./final11.pstex_t}}} \caption{}
\label{fig4}
\end{figure}

Since $R$ is arbitrary, we get that $u_\infty(x,y):=\varphi(y) \leq
u(x,y)$ for all $(x,y) \in H$. Moreover, the normal derivative of
the functions $u$ and $u_\infty$ is the same at the origin, and by
the maximum principle we get
\[
u=u_\infty\,.
\]
This shows that $\Omega=H$. Therefore we just need to show the claim.

\medskip

{\it Proof of the claim:} Fix $R>0$ and fix $B$ a ball of radius
$2R$ tangent to a certain point $q_n$, with sufficiently large
$n$. By Lemma \ref{np}, inward normal half-line $N(q_n)$ starting
at $q_n$ does not intersect $\partial \Omega$. Moreover, at a
certain point it reaches the half-plane $H$.

We move the center of $B$ along $N(q_n)$ till it reaches $H$. We
first show that during that motion, $B$ cannot intersect $\partial
\Omega$ at both sides of $N(q_n)$. Indeed, denote by $p_1$ and
$p_2$ such intersection points. Then, the ball
$B_{|p_1-p_2|+1}(p_1)$ would give a contradiction with Lemma
\ref{otro}.

Therefore, when we move the center of $B$ along $N(q_n)$, it
eventually intersects $\partial \Omega$ just from one side.
Therefore we can move a ball of radius $R$ up to the half-plane
$H$ through the other side. From there, we can easily translate it
to reach the position of $B_R(q)$.
\proofend

\enddocument